\documentclass{article}
\usepackage[utf8]{inputenc}

\title{Another Enumeration of Caterpillar Trees}
\author{Jacob Crabtree}
\date{October 27, 2018}
\usepackage{graphicx}
\usepackage{float}
\usepackage{amsmath}
\usepackage{amsfonts}
\usepackage{mathtools}

\begin{document}

\maketitle

\section*{}
\underline{\textbf{Abstract:}} A caterpillar tree is a connected, acyclic, graph in which all vertices are either a member of a central path, or joined to that central path by a single edge. In other words, caterpillar trees are the class of trees which become path graphs after removing all leaves. In 1973, F. Harary and A.J. Schwenk provided two proofs found in [1] which show that the number of non-isomorphic caterpillars with N vertices is given by the formula $2^{N-4} + \ 2^{ \lfloor \frac{N - 4}{2}\rfloor}$, where $\lfloor \  \rfloor$ denotes the floor function. The first proof follows from a special case of an application of P\'{o}lya's Enumeration theorem on graphs with integer-weighted vertices. The second proof proceeds through an appropriate edge labelling of the caterpillars. The proof presented here owes much of its insight to the first two, but has the benefit of utilizing a natural labelling for the caterpillars. We will proceed by labelling the vertices of the caterpillars with integer-weights, followed by an application of the orbit-counting theorem.

\section*{Weighted Path Graphs and Spines}

One advantage of the class of caterpillar trees is their succinct representation as a weighted path graph. A weighted path graph can be obtained from a caterpillar by first assigning an integer label to each vertex in its 'central path' representing the number of leaves adjacent to it. Then, all of the leaves can be removed from the graph, with the exception of two designated endpoints which are labeled with zeros and treated as members of the central path. This process is depicted in figure 1.

\begin{figure}[H]
\centering
\includegraphics[scale=0.6]{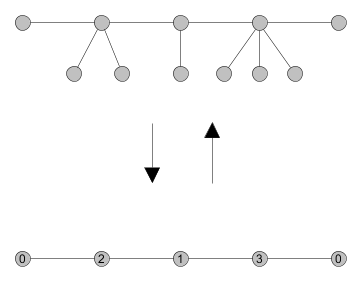}
\caption{Caterpillar labeled as a weighted graph}
\end{figure}

At first glance, it may seem quite unnatural to include the two zeros in our labeling. However, the necessity of these zeros will become clear after some additional considerations. We begin with a definition.\bigbreak
\underline{\textbf{\textit{Definition:}}} \textit{Let $T \ \text{=} \ \{V, E\}$ be a labelled caterpillar tree with vertex set\newline $V \ \text{=} \ \{v_1,...,v_n\}$ and edge set $E \ \text{=} \ \{e_1,...,e_m\}.$ A \textbf{derived path}, $D \ \text{=} \ \{V',E'\}$, of the caterpillar $T$ is a subtree, such that $E'$ is a path in $T$, and for each path $P \subset E$, we have $|P| \leq |E'|.$}\bigbreak

Quite simply, a derived path is a maximal path in $T$, and different derived paths are simply different designations of endpoints. However, as mentioned, when we consider the weighted representation of the tree, the vertices which we choose to designate as endpoints are always assigned the label $0.$ Moreover, since the derived paths may differ only by endpoints, the remaining integer labels assigned to non-endpoints will be the same in any derived path for a given caterpillar tree. With this in mind, we define the \textbf{\textit{spine}} of a caterpillar as follows:\bigbreak

\underline{\textbf{\textit{Definition:}}} \textit{Let $T \ \text{=} \ \{V,E\}$ be the caterpillar defined above, and let $D$ be a derived path of $T$ with vertex set $V' \subset V$ and edge set $E' \subset E$. For simplicity, we will assume that $V' \ \text{=} \ \{v_1,v_2,...v_k\},$ where $v_1$ and $v_k$ are the two endpoints. We will also assume that $D$ has been drawn in a horizontal line with $v_1$ on the far left, $v_k$ on the far right, and the remaining vertices labeled from left to right. Let $N(v_i)$ denote the neighborhood of the vertex $v_i.$ Then, the \textbf{left-right oriented spine} of the caterpillar $T$ is a vector $\overline{s} \ \text{=} \ (s_2,s_3,...,s_{k-1})$ of dimension $k-2$ such that the component $s_i \ \text{=}  \ |N(v_i) - V'|.$ In other words, $s_i$  represents the number of vertices adjacent to vertex $v_i$ which are not themselves members of the derived path.}\bigbreak

For example, the caterpillar from figure 1 obtains the spine $(2,1,3).$ This procedure allows us to obtain a left-right oriented spine from a given caterpillar tree. In the other direction, we can obtain a weighted tree from a given spine $\overline{s} \ \text{=} \ (s_1,...,s_k)$ by drawing a path graph of length $k,$ and labeling the vertices from left to right with the corresponding entries in the spine. We then add an additional vertex onto both ends of the path and label these two new vertices with a $0$. This weighted tree now corresponds to some caterpillar tree, $C$, which may be obtained by attaching the labeled number of flags at each vertex. We will refer to this reverse process by saying that the left-right oriented caterpillar $C$ is \textbf{\textit{induced}} by the left-right oriented spine $\overline{s}.$ In all subsequent discussion we will agree to always use left-right oriented spines to draw left-right oriented caterpillars so that we may refer to \textit{the} caterpillar and \textit{the} spine without ambiguity.

As mentioned above, one may be tempted to define the spine by removing \textit{all} leaves from a particular caterpillar instead of using any zeros or considering derived paths. For example, this 'naive scheme' would assign the caterpillar $C_1$ from figure 1 the spine $\overline{s_1} \ \text{=} \ (3,1,4).$ However, the naive scheme would allow spine vectors which have the same length to induce caterpillars which have different maximal path lengths. Consider for example the caterpillar $C_2$ induced by the spine $\overline{s_2} \ \text{=} \ (0,4,4).$ Although $\overline{s_1}$ and $\overline{s_2}$ are the same length as vectors, the caterpillars $C_1$ and $C_2$ induced under this naive scheme would have maximal path lengths of 5 and 4 respectively. Our definition for the spine does not allow for such confusion. By considering the endpoints separately from the spine, we ensure that all spines of length $k$ induce caterpillars which have maximal path length $k+2.$ We will now proceed with another definition.\bigbreak

\underline{\textbf{\textit{Definition:}}} \textit{The \textbf{spine class}, which we denote $\mathcal{S}_{N_k}$, is the set of spines of length $k$ which induce caterpillar trees with $N$ vertices (and maximal path length $k+2$). That is,} $$\mathcal{S}_{N_k} \ \text{=} \ \bigg\{(s_1,...,s_k) \ \Big| \ \bigg(\sum_{i = 1}^{k} s_i \bigg) \ \text{=} \ N - k  - 2 \bigg\}.$$\bigbreak

Before proceeding, we make one last remark about the relationship between spines and caterpillars. Given a caterpillar tree, we can assign a unique spine to it by following the procedure outlined in the definition of a spine. However, given two spines, $\overline{s_1}, \overline{s_2}$ such that $\overline{s_1} \not = \overline{s_2}$, the respective induced caterpillars, $C_1, C_2,$ are not necessarily distinct. For example, let $\overline{s_1} \ \text{=} \ (1,1,2)$, $\overline{s_2} \ \text{=} \ (2,1,1)$. The spines and the corresponding caterpillars are depicted in figure 2. 

\begin{figure}[H]
\centering
\includegraphics[scale=0.20]{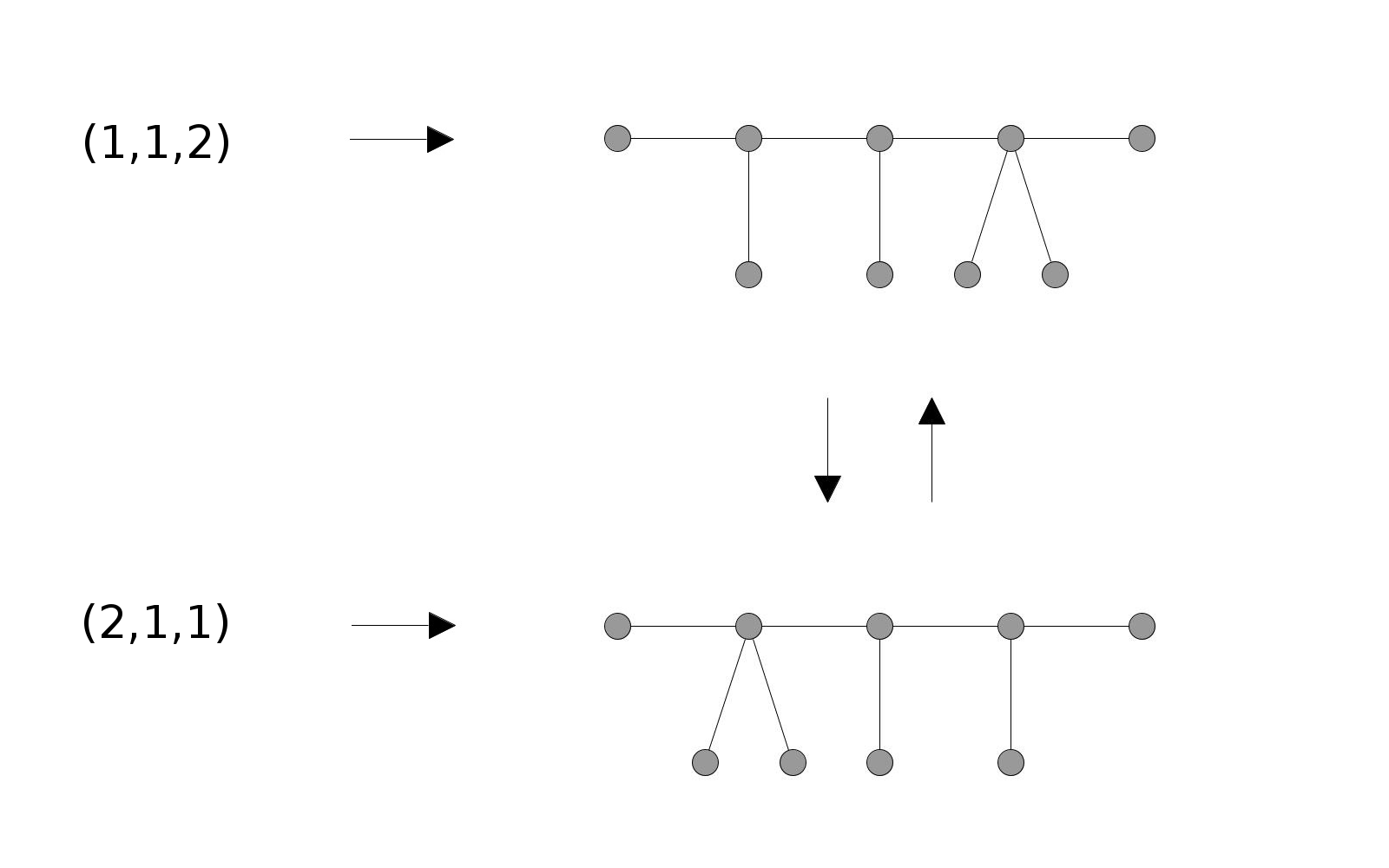}
\caption{Two distinct spines which induce isomorphic caterpillar trees}
\end{figure}

Although the two spines from figure 2 are distinct as vectors, the two induced caterpillars are clearly isomorphic and are related by reflection over a vertical line centered on the middle vertex. Obviously this is only an artifact of our choice of left-right orientation for the spine vectors and graph drawings. We can frame this observation more clearly in the language of group theory. Let $\mathcal{S}_{N_k}$ denote the spine class as defined above, and let $\mathbb{Z}_2$ denote the cyclic group of order two. Then, we may define a group action of $\mathbb{Z}_2$ on $\mathcal{S}_{N_k}$ in a natural way. Let $\overline{s} \ \text{=} \ (s_1,s_2,...,s_{k-1},s_k) \in \mathcal{S}_{N_k} $ and let $\alpha$ denote the non-identity element of $\mathbb{Z}_2$. Then, we define the action $\alpha * \overline{s} \ \text{=} \ \overline{s'}$ by $s'_i \ \text{=} \ s_{k-i+1}$, for $1 \leq i \leq k.$ It is clear that two caterpillars $C_1, C_2$, with respective spines $\overline{s_1}, \overline{s_2}$ are isomorphic if and only if $s_1 \ \text{=} \ g * \overline{s_2},$ for some $g \in \mathbb{Z}_2.$ Therefore, in order to count the number of caterpillars with $N$ vertices, and spine length $k$, we can instead count the orbits of the spine class $\mathcal{S}_{N_k}$ under the action of $\mathbb{Z}_2$. To this end, we record here the well known orbit-counting theorem. A proof can be found in [2, p. 196-199].\bigbreak 

\underline{\textbf{\textit{Theorem:}}} \textit{Let $X$ be any set and let $G$ be a finite group that acts on $X$. We will denote the set of orbits under this action by $X/G$, and the invariant set of an element of a group by $\big(X\big)^g.$ That is, for some element $g \in G,$ we define $\big(X\big)^g \ \text{=} \ \{x \in X \ | \ g*x \text{=} \ x\}.$ The orbit-counting theorem states the following:}
$$|X/G| \ \text{=} \ \frac{1}{|G|} \sum_{g \in G} |\big(X\big)^g|$$\bigbreak

\section*{Three Bijections}

We have described how a particular set of vectors corresponds to of a particular class of trees. Our new question becomes, "how many distinct $k$ dimensional vectors have the property that the sum of the $k$ components is equal to a fixed constant?". Thankfully, this revised question has been well investigated. In fact, if we restrict ourselves to the case of vectors with strictly positive entries, then the number of distinct vectors is given by the famed Stirling numbers of the second kind. One can prove this theorem, as well as a number of others, through a wonderful style of argument known as 'Stars and Bars Combinatorics'. We record one theorem which will be of use to us here. A proof can be found in [3, p. 12-13].\bigbreak

\underline{\textbf{\textit{Theorem:}}} \textit{For any pair of positive integers n and k, the number of k-tuples of non-negative integers whose sum is n is equal to the binomial coefficient ${n + k - 1}\choose{k - 1} $}

By substituting $n \ \text{=} \ N - k - 2$, we obtain an equation for the order of a given spine class:
\begin{align}
|\mathcal{S}_{N_k}| \ \text{=} \  \binom{N - 3}{k - 1}
\end{align}

Obviously the order of a given spine class also gives the order of the invariant set for the identity element, $\big(\mathcal{S}_{N_k}\big)^e$. In order to make use of the orbit-counting theorem, what remains is to count the invariant set of the non-identity element, $\big(\mathcal{S}_{N_k}\big)^\alpha$. In other words, we must count the number of caterpillars which are their own reflection. Harary and Schwenk refer to members of this invariant set as the \textit{symmetric caterpillars} [1], and made the critical observation required to count them. Namely, any symmetric caterpillar can be constructed from two copies of some smaller caterpillar. There are four cases we must consider which correspond to the four possible combinations of parity for $N$ and $k$. As we will see shortly, there is a correspondence between three of the four cases, and the final case is impossible. We will begin with an illustration of the correspondence.

\begin{figure}[H]
\centering
\includegraphics[scale=0.175]{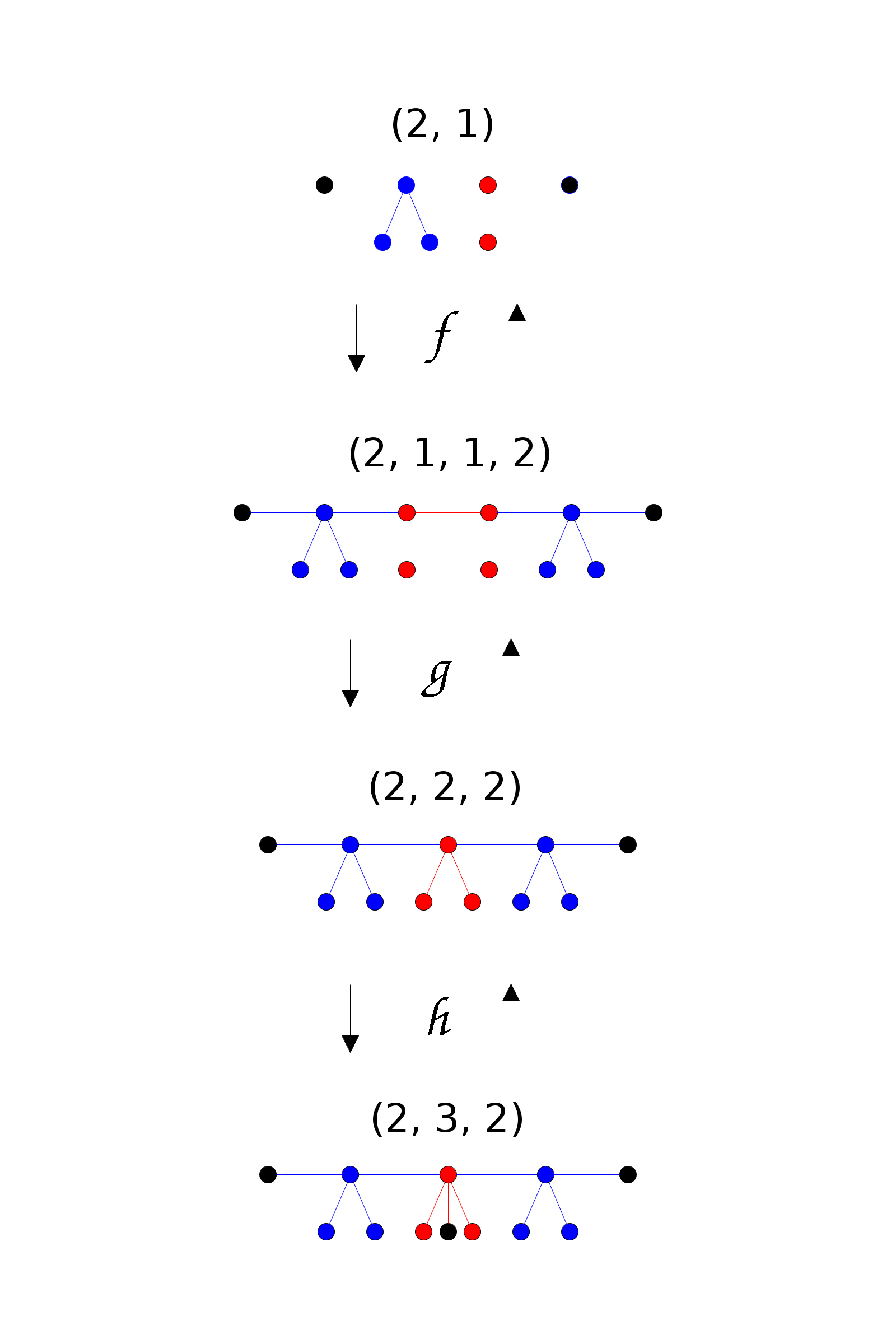}
\caption{Spine class bijections for three sets of symmetric caterpillars.}
\end{figure}

As figure 3 suggests, there is a bijective correspondence between the three sets of symmetric caterpillars with a certain class of 'half-spines'. Starting from the top of figure three, we see that the first two caterpillars are related by taking a copy of the 'half-spine', reflecting it, and then gluing the two copies together. The second and third graphs are related by expanding or contracting the center edge. The third and fourth graphs are related by the addition or deletion of a single flag from the center vertex.

While a careful examination of figure 3 is likely the best way to understand the trick at the core of our main proof, we can also construct the bijections $f,g,$ and $h$ more explicitly. In the following proofs $N$ and $k$ will always denote an even number. When there is a need to reference an odd $N$ or odd $k$ we will do so by using $N \pm 1$ and $k \pm 1$ as necessary.

Let $\overline{x} \ \text{=} \ (x_1,...,x_{k}),$ and $\overline{y} \ \text{=} \ (y_1,...,y_k)$ be two spine vectors in $\big(\mathcal{S}_{N_k}\big)^\alpha$. Define the function $f:\big(\mathcal{S}_{N_k}\big)^\alpha \rightarrow \mathcal{S}_{(\frac{N}{2} + 1)_{\frac{k}{2}}}$ by $f(\overline{x}) \ \text{=} \ (x_1, ... , x_{\frac{k}{2}}).$  The claim is that $f$ is the required bijection.\bigbreak

\underline{\textbf{\textit{Proof:}}}  To see that $f$ is one-to-one, suppose that $f(\overline{x}) \ \text{=} \ f(\overline{y}).$ Then, for $i \leq \frac{k}{2},$ we have $x_i \ \text{=} \ y_i.$ However, by the definition of $\big(\mathcal{S}_{N_k}\big)^\alpha,$ we also have that $x_i \ \text{=} \ x_{k-i+1},$ and $y_i  \ \text{=} \ y_{k-i+1},$ which implies that $x_i  \ \text{=} \ y_i$ for $\frac{k}{2} \leq i \leq k$ as well. Hence, we see that $x_i  \ \text{=} \ y_i$ for all $1 \leq i \leq k$ so that $\overline{x} \ \text{=} \ \overline{y}.$ Next we will show that $f$ is onto. Suppose $\overline{z} \ \text{=} \ (z_1, z_2, ...,z_{\frac{k}{2}}) \in \mathcal{S}_{(\frac{N}{2} + 1)_{\frac{k}{2}}}.$ We must produce a vector $\overline{x} \in \big(\mathcal{S}_{N_k}\big)^\alpha$ such that $f(\overline{x}) \ \text{=} \ \overline{z}.$ To this end, we define the candidate vector $\overline{x} \ \text{=} \ (z_1, z_2 ... ,z_{\frac{k}{2}},z_{\frac{k}{2}}, ... ,z_2 , z_1).$ First, we see that $x$ has length $k$ and that $\alpha * \overline{x} = \overline{x}$ by its construction. Next, we see that $$ \sum_{i = 1}^{k} x_i \ \text{=} \ 2\sum_{i = 1}^{\frac{k}{2}} z_i \ \text{=} \ 2\Bigg[ \bigg(\frac{N}{2} + 1\bigg) - \bigg(\frac{k}{2}\bigg) - 2 \Bigg] \ \text{=} \ N - k -2. $$ Therefore, we have ensured that $\overline{x} \in \big(\mathcal{S}_{N_k}\big)^\alpha,$ so that $\overline{x}$ serves as the sought after vector. Hence, $f$ is also onto, which completes the proof.

\hfill $\square.$\bigbreak

The proof of the next two cases is very similar. Let $\overline{x} \ \text{=} \ (x_1,...,x_{k}),$ and $\overline{y} \ \text{=} \ (y_1,...,y_k)$ be two spine vectors in $\big(\mathcal{S}_{N_k}\big)^\alpha$. Define the function $g:\big(\mathcal{S}_{N_k}\big)^\alpha \rightarrow \big(\mathcal{S}_{(N-1)_{k-1}}\big)^\alpha$ by
\begin{align*}
g(\overline{x}) \ &\text{=} \ g((x_1,x_2, ... ,x_{\frac{k}{2} - 1}, x_{\frac{k}{2}},x_{\frac{k}{2} + 1}, x_{\frac{k}{2} + 2}, ... , x_{k-1}, x_k))\\
\ &\text{=} \ (x_1,x_2,...,x_{\frac{k}{2} - 1}, x_{\frac{k}{2}} + x_{\frac{k}{2} + 1}, x_{\frac{k}{2} + 2} ..., x_{k-1}, x_{k})
\end{align*}
We wish to show that $g$ is a bijection.\bigbreak

\underline{\textbf{\textit{Proof:}}} First we show that $g$ is one-to-one. The proof is essentially the same as the proof given for $f$ above. Let $g(\overline{x}) \ \text{=} \ g(\overline{y}).$ Then, we can see immediately that $x_i = y_i,$ for $1 \leq i \leq \frac{k}{2} - 1,$ and $\frac{k}{2} + 2 \leq i \leq k.$ We can also see that $x_{\frac{k}{2}} + x_{\frac{k}{2} + 1} \ \text{=} \ y_{\frac{k}{2}} + y_{\frac{k}{2} + 1}
.$ However, since $\overline{x}$ and $\overline{y}$ are symmetric, $x_{\frac{k}{2}} \ \text{=} \ x_{\frac{k}{2} + 1}$ and $y_{\frac{k}{2}} \ \text{=} \ y_{\frac{k}{2} + 1}.$ Therefore, we may simplify to $2x_{\frac{k}{2}} \ \text{=} \ 2y_{\frac{k}{2}},$ which implies $x_{\frac{k}{2}} \ \text{=} \ x_{\frac{k}{2} + 1} \ \text{=} \ y_{\frac{k}{2} + 1} \ \text{=} \ y_{\frac{k}{2}}.$ Hence, $\overline{x} \ \text{=} \ \overline{y}$ so that $g$ is one-to-one. Next, we will show that $g$ is onto. To begin, we will show first that each $\overline{z} \in \big(\mathcal{S}_{(N-1)_{k-1}}\big)^\alpha$ can be written in the form $(z_1,z_2,...,z_{\frac{k}{2} - 1}, 2q, z_{\frac{k}{2} + 1} ..., z_{k-2}, z_{k-1}),$ for some value of $q.$ Suppose conversely that for some $\overline{z}$ we have $$\overline{z} \ \text{=} \ (z_1,z_2,...,z_{\frac{k}{2} - 1}, 2q + 1, z_{\frac{k}{2} + 1} ..., z_{k-2}, z_{k-1}).$$ 
Since $\overline{z}$ is symmetric, $z_{i} = z_{k-i-1}$. Taking the sum of the components, we see that $$\sum_{i = 1}^{k-1} z_i\ \text{=} \ (2q + 1) + 2\sum_{i = 1}^{\frac{k}{2} - 1}z_i, $$
which is odd. However, from the definition of the spine class $\mathcal{S}_{(N-1)_{k-1}},$ we also have: $$ \sum_{i = 1}^{k-1} z_i \ \text{=} (N-1) - (k-1) - 2, $$
which is even. This contradiction shows that each $\overline{z} \in \big(\mathcal{S}_{(N-1)_{k-1}}\big)^\alpha$ must have the form $(z_1,z_2,...,z_{\frac{k}{2} - 1}, 2q, z_{\frac{k}{2} + 1} ..., z_{k-2}, z_{k-1}),$ for some value of $q.$ To produce the vector $\overline{x} \in \big(\mathcal{S}_{N_k}\big)^\alpha$ such that $g(\overline{x}) \ \text{=} \ \overline{z},$ we simply set $$\overline{x} = (z_1,z_2,...,z_{\frac{k}{2} - 1},q,q,z_{\frac{k}{2} + 1},...,z_{k-2},z_{k-1}).$$ Since the vector $\overline{z}$ has length $k-1$ it is clear that The vector $\overline{x}$ has length $k$. Furthermore, we see that  $$\sum_{i = 1}^{k} x_i \ \text{=} \ \sum_{i = 1}^{k-1} z_i \ \text{=} \ c.$$ Since $\overline{y} \in \mathcal{S}_{(N-1)_{k-1}},$ we have that $c = (N-1) - (k-1) - 2$. We also have that $c = (N_x) - k - 2$ for some $N_x$ which corresponds to the spine class which $\overline{x}$ belongs to. Clearly this gives $N_x = N$ so that $\overline{x} \in \mathcal{S}_{N_{k}}.$ It is also easy to see that by construction, $\alpha * \overline{x} = \overline{x}$ so that $x \in \big(\mathcal{S}_{N_{k}}\big)^\alpha \ \text{and} \ g(\overline{x}) = \overline{y}.$ Therefore, we have shown that $g$ is also onto.

\hfill $\square.$\bigbreak

Finally, we note that the function $h:\big(\mathcal{S}_{(N-1)_{k-1}}\big)^\alpha \rightarrow \big(\mathcal{S}_{N_{k-1}}\big)^\alpha$ given by 
$$h(x_i) \ \text{=} \
\begin{cases} 
      x_{\frac{k}{2}} + 1, &  i \ \text{=} \ \frac{k}{2}\\
      x_i, &  \text{else} \\
   \end{cases}
$$
and the inverse function  $h^{-1}:\big(\mathcal{S}_{N_{k-1}}\big)^\alpha \rightarrow \big(\mathcal{S}_{(N-1)_{k-1}}\big)^\alpha $ given by 
$$h^{-1}(z_i) \ \text{=} \
\begin{cases} 
      z_{\frac{k}{2}} - 1, &  i \ \text{=} \ \frac{k}{2}\\
      z_i, &  \text{else} \\
   \end{cases}
$$
defines the final bijection. The proof is essentially the same as those for $f$ and $g$.\newpage 

These three bijections yield the following result:
\begin{align}
    \Big|\big(\mathcal{S}_{N_k}\big)^\alpha \Big| \ \text{=} \ \Big|\big(\mathcal{S}_{N_{k-1}}\big)^\alpha \Big| \  \text{=} \ \Big|\big(\mathcal{S}_{{(N-1)}_{k-1}}\big)^\alpha \Big| \ \text{=} \ \Big|\mathcal{S}_{(\frac{N}{2} + 1)_{\frac{k}{2}}} \Big|
\end{align}
Combining (1) and (2) we obtain:
\begin{align}
    \Big|\big(\mathcal{S}_{N_k}\big)^\alpha \Big| \text{=} \Big|\big(\mathcal{S}_{N_{k-1}}\big)^\alpha \Big| \text{=} \Big|\big(\mathcal{S}_{{(N-1)}_{k-1}}\big)^\alpha \Big| \text{=} \Big|\mathcal{S}_{(\frac{N}{2} + 1)_{\frac{k}{2}}} \Big| \text{=} \binom{\frac{N}{2} - 2}{\frac{k}{2} - 1}
\end{align} 
The formula (3) accounts for three of the four cases of possible parities for $N$ and $k$. The remaining case is for $N-1$, $k$ (i.e. odd number of vertices, even spine length). In this case, a simple argument shows that $\big(\mathcal{S}_{{(N-1)}_k}\big)^\alpha$ \ \text{=} \ $\emptyset$.\bigbreak

\underline{\textbf{\textit{Proof:}}} Suppose there were some $\overline{s} \ \text{=} \ (s_1,...,s_k) \in \big(\mathcal{S}_{{(N-1)}_k}\big)^\alpha$. By definition of the invariant set $\big(\mathcal{S}_{{(N-1)}_k}\big)^\alpha,$ we have $s_i \ \text{=} \ s_{k-i+1}$, for $1 \leq i \leq k.$ But then $$ \sum_{i = 1}^{k} s_i \ \text{=} \ 2\sum_{i = 1}^{k/2} s_i, $$ which is even. However, by definition, of the spine class $\mathcal{S}_{{(N-1)}_k}$, we also have $$ \sum_{i = 1}^{k} s_i \ \text{=} \ (N-1) - k - 2 $$ which is odd.

\hfill $\square.$\bigbreak 

This contradiction shows that in the case of $N-1$, $k$, we must have

\begin{align}
\big(\mathcal{S}_{{(N-1)}_k}\big)^\alpha \ \text{=} \ \emptyset.
\end{align}\bigbreak

\section*{Proof of Main Theorem:}
With these observations in hand, the stage has been set and we are ready to prove our main theorem. Let $\mathcal{C}_N$ denote the set of caterpillars with $N$ vertices. In order to evaluate $|\mathcal{C}_N|,$ we will evaluate $|\mathcal{S}_{N_k} / \mathbb{Z}_2|$ for each legal spine length $k,$ and then take the summation over all such $k$. The definition of a caterpillar requires that the minimum spine length be 1, and the definition of the spine requires that the maximum spine length be $N-2$. Therefore, these will serve as the bounds for our summation. We will begin with the case where $N$ is even: 

\begin{equation*} 
\begin{split}
|\mathcal{C}_N| \ &\text{=} \  \sum_{i = 1}^{N-2} |\mathcal{S}_{N_i} / \mathbb{Z}_2| \\
\ &\text{=} \ \sum_{i = 1}^{N-2} \frac{1}{|\mathbb{Z}_2|} \Bigg( \big|\big(\mathcal{S}_{N_i}\big)^e\big| + \big|\big(\mathcal{S}_{N_i}\big)^\alpha\big| \Bigg) \\
\ &\text{=} \ \frac{1}{2} \Bigg( \sum_{i = 1}^{N-2} \big|\mathcal{S}_{N_i}\big| + \sum_{i = 1}^{N-2} \big|\big(\mathcal{S}_{N_i}\big)^\alpha\big| \Bigg)  \\
\ &\text{=} \ \frac{1}{2} \Bigg( \sum_{i = 1}^{N-2} \big|\mathcal{S}_{N_i}\big| + \sum_{i = 1,3,5,...}^{N-2} \big|\big(\mathcal{S}_{N_i}\big)^\alpha\big| + \sum_{i = 2,4,6,...}^{N-2} \big|\big(\mathcal{S}_{N_i}\big)^\alpha\big| \Bigg)  \\
\ &\text{=} \ \frac{1}{2} \Bigg( \sum_{i = 1}^{N-2} \big|\mathcal{S}_{N_i}\big| + 2\sum_{i = 2,4,6,...}^{N-2} \big|\big(\mathcal{S}_{N_i}\big)^\alpha\big| \Bigg)  \\
\ &\text{=} \ \frac{1}{2} \Bigg( \sum_{i = 1}^{N-2} \big|\mathcal{S}_{N_i}\big| + 2\sum_{i = 2,4,6,...}^{N-2} \bigg|\mathcal{S}_{(\frac{N}{2} + 1)_{\frac{i }{2}}}\bigg| \Bigg)  \\
\ &\text{=} \ \frac{1}{2} \Bigg( \sum_{i = 1}^{N-2} \big|\mathcal{S}_{N_i}\big| + 2\sum_{i = 1}^{\frac{N}{2} - 1} \bigg|\mathcal{S}_{(\frac{N}{2} + 1)_{i}}\bigg| \Bigg) \\
\ &\text{=} \ \frac{1}{2} \Bigg( \sum_{i = 1}^{N-2} \binom{N - 3}{i - 1} + 2\sum_{i = 1}^{\frac{N}{2} - 1} \binom{\frac{N}{2} - 2}{i - 1} \Bigg) \\
\ &\text{=} \ \frac{1}{2} \Bigg( \sum_{i = 0}^{N-3} \binom{N - 3}{i} + 2\sum_{i = 0}^{\frac{N}{2} - 2} \binom{\frac{N}{2} - 2}{i} \Bigg) \\
\ &\text{=} \ \frac{1}{2} \bigg( 2^{N-3} + 2*2^{\frac{N}{2} - 2} \bigg) \\
\ &\text{=} \ 2^{N-4} + 2^{\frac{N-4}{2}} \\
\end{split}
\end{equation*}

\hfill $\square.$\bigbreak

Here we have made use of the orbit-counting theorem from lines one-two, and made use of theorem (3) from lines four-five, five-six, and seven-eight. From lines nine-ten we made use of the well known formula for the sum of binomial coefficients. 

We now turn our attention to the case where $N$ is odd. The order of the invariant sets $\big(\mathcal{S}_{N_i}\big)^e$ will be the same as before. By making use of our counting theorem (3), we can make a simple adjustment to the formula for the order of the invariant set $\big(\mathcal{S}_{N_i}\big)^\alpha.$\newpage 
Suppose $N$ and $k$ are both odd. Then $N+1$ and $k+1$ are both even. By applying (3) we obtain:
\begin{align*}
    \Big|\big(\mathcal{S}_{{(N+1)}_{(k+1)}}\big)^\alpha \Big| \ \text{=}
    \Big|\big(\mathcal{S}_{{N}_{k}}\big)^\alpha \Big| \ \text{=} \ \Big|\mathcal{S}_{(\frac{N+1}{2} + 1)_{\frac{k+1}{2}}} \Big| \ \text{=} \ \binom{\frac{N+1}{2} - 2}{\frac{k+1}{2} - 1}.
\end{align*}
If $N$ is odd but $k$ is even, then by (4) we have $\Big|\big(\mathcal{S}_{{(N)}_{k}}\big)^\alpha \Big| \ \text{=} \ 0.$ Making these adjustments, we obtain the formula for case of odd $N$:

\begin{equation*} 
\begin{split}
|\mathcal{C}_N| \ &\text{=} \  \sum_{i = 1}^{N-2} |\mathcal{S}_{N_i} / \mathbb{Z}_2| \\
\ &\text{=} \ \sum_{i = 1}^{N-2} \frac{1}{|\mathbb{Z}_2|} \Bigg( \big|\big(\mathcal{S}_{N_i}\big)^e\big| + \big|\big(\mathcal{S}_{N_i}\big)^\alpha\big| \Bigg) \\
\ &\text{=} \ \frac{1}{2} \Bigg( \sum_{i = 1}^{N-2} \big|\mathcal{S}_{N_i}\big| + \sum_{i = 1}^{N-2} \big|\big(\mathcal{S}_{N_i}\big)^\alpha\big| \Bigg)  \\
\ &\text{=} \ \frac{1}{2} \Bigg( \sum_{i = 1}^{N-2} \big|\mathcal{S}_{N_i}\big| + \sum_{i = 1,3,5,...}^{N-2} \big|\big(\mathcal{S}_{N_i}\big)^\alpha\big| + \sum_{i = 2,4,6,...}^{N-2} \big|\big(\mathcal{S}_{N_i}\big)^\alpha\big|  \Bigg)\\
\ &\text{=} \ \frac{1}{2} \Bigg( \sum_{i = 1}^{N-2} \big|\mathcal{S}_{N_i}\big| + \sum_{i = 1,3,5,...}^{(N+1)-2} \bigg|\mathcal{S}_{(\frac{N+1}{2} + 1)_{\frac{i + 1}{2}}}\bigg| + 0 \Bigg)  \\
\ &\text{=} \ \frac{1}{2} \Bigg( \sum_{i = 1}^{N-2} \binom{N - 3}{i - 1} + \sum_{i = 1,3,5,...}^{(N+1)-2} \binom{\frac{N+1}{2} - 2}{\frac{i + 1}{2} - 1} \Bigg) \\
\ &\text{=} \ \frac{1}{2} \Bigg( \sum_{i = 1}^{N-2} \binom{N - 3}{i - 1} + \sum_{i = 1}^{\frac{N+1}{2} - 1} \binom{\frac{N+1}{2} - 2}{i - 1} \Bigg) \\
\ &\text{=} \ \frac{1}{2} \Bigg( \sum_{i = 0}^{N-3} \binom{N - 3}{i} + \sum_{i = 0}^{\frac{N+1}{2} - 2} \binom{\frac{N+1}{2} - 2}{i} \Bigg) \\
\ &\text{=} \ \frac{1}{2} \bigg( 2^{N-3} + 2^{\frac{N+1}{2} - 2} \bigg) \\
\ &\text{=} \ 2^{N-4} + 2^{\frac{N-5}{2}} \\
\end{split}
\end{equation*}

\hfill $\square.$ \bigbreak

Here we applied our counting theorems (3) and (4) from lines four-five, and five-six. Finally, we can use the floor function, $\lfloor x \rfloor$, to combine these two results into a succinct formula which gives the number of non-isomorphic caterpillar trees for any prescribed number of vertices, $N$: $$|\mathcal{C}_N| \ \text{=} \  2^{N-4} + 2^{\big\lfloor\frac{(N-4)}{2}\big\rfloor}.$$ \newpage

\section*{}
\textbf{References:}\bigbreak
[1] Harary, F., \& Schwenk, A. J. (1973). The number of caterpillars. \textit{Discrete Mathematics}, 6(4), 359-365. doi:10.1016/0012-365x(73)90067-8\bigbreak
[2] Harris, J. M., Hirst, J. L., \& Mossinghoff, M. J. (2008). \textit{Combinatorics and graph theory}. New York: Springer.\bigbreak
[3] Ross, S. M. (2010). \textit{A first course in probability.} Boston: Pearson. 

\end{document}